\def\pitilde{\tilde\pi(\zeta|z)}
\def\ptilde{\tilde p(z|\zeta)}

\def\pipt{{\pi_{\mathrm{pt}}}}
\def\pipr{{\pi_{\mathrm{prob}}}}
\def\onehalf{{\textstyle{1\over2}}}
\documentclass[final] {aipproc}

\layoutstyle{6x9}

\begin{document}

\title{The Marginalization Paradox \\
and the Formal Bayes' Law}

\classification{02.50.-r, 02.50.Tt, 02.70.Rr}
\keywords      {Marginalization paradox, Objective Bayes, Logical Bayes,
Improper Priors, Maximum Entropy}

\author{Timothy C. Wallstrom}{
  address={Theoretical Division, MS B213 \\
    Los Alamos National Laboratory \\
    Los Alamos, NM 87545 \\
\texttt{tcw@lanl.gov}}
}

\begin{abstract}
  It has recently been shown that the marginalization paradox~(MP) can
  be resolved by interpreting improper inferences as probability
  limits.  The key to the resolution is that probability limits need
  not satisfy the formal Bayes' law, which is used in the MP to deduce
  an inconsistency.  In this paper, I explore the differences between
  probability limits and the more familiar pointwise limits, which do
  imply the formal Bayes' law, and show how these differences underlie
  some key differences in the  interpretation of the MP.
\end{abstract}

\maketitle

\section{Introduction}
\label{sec:intro}

The marginalization paradox~(MP) is an apparent inconsistency in
Bayesian inference that can arise from the use of improper priors. It
was discovered in 1972 by Dawid and Stone~\cite{SD72}; together with
Zidek, they published a comprehensive analysis in
1973~\cite{DSZ73}. We follow Jaynes in referring to these authors as
``DSZ.''

The MP arises in problems with a particular structure, where there are
two different ways of computing the same marginal posterior. When
improper priors are used, the two results are usually incompatible; an
inconsistency cannot arise with proper priors.  The improper
inferences are computed as ``formal posteriors'' using the usual
Bayes' law, $\pi(\theta|x)\propto p(x|\theta)\,\pi(\theta)$, with an
improper prior.

The MP is important to objective Bayes because the ``noninformative''
priors required by the theory are typically improper. By suggesting
that improper priors cannot be used consistently, the MP raises
questions as to whether noninformative priors exist, and thus as to
whether the objective Bayesian approach is tenable. We use the term
objective Bayes to describe the approach of Harold Jeffreys and Edwin
Jaynes, in which a certain state of information is represented by a
unique prior. (Their approach is also known as ``logical Bayes''; the
term objective Bayes is now often used for a different but related
approach in which the prior may depend on the
estimand~\cite{Bernardo79}; for a discussion of different approaches,
cf.~\cite{KW96}.)

Edwin Jaynes strongly contested the view that the MP represented a
true inconsistency, and engaged DSZ in a spirited debate;
cf.~\cite{Jaynes03} and references therein. The battle lines were
drawn in the 1970s, and have changed little since. Neither side
managed to convince the other, and the absence of new results has
caused the paradox to be largely set aside, even though it has never
been completely understood.

Recently, I have shown~\cite{Wallstrom07} that the MP can be resolved
if probability limits, rather than formal posteriors, are used to
define improper inferences. The purpose of this paper is to reconsider
the differences between Jaynes and DSZ in the light of this new
result. One might assume, in reviewing their debates, that Jaynes and
DSZ were separated by an unbridgeable chasm.  It is shown, by
contrast, that the differences between Jaynes and DSZ hinge on a
single assumption, which might at first appear as a mere technicality.
This assumption is the type of limit used to define the improper
inference. We then make the case that the limit process used to
support the use of formal posteriors is unsound. Our analysis relies
heavily on the ideas of Mervyn Stone, who first introduced the notion
of probability limit, and has long maintained that the paradoxes
associated with formal posteriors are due to the inadequacy of the
pointwise limit.

We also discuss the impact of these recent findings on the prospects
for objective Bayes. One might have hoped that a resolution of the
paradox might open the door to a refined theory of improper inference,
using probability limits instead of formal posteriors, which would
provide a consistent foundation for objective Bayes. Unfortunately,
there is strong evidence that probability limits rarely exist.
Although we can identify some problems that we can now solve which
previously led to inconsistency, it appears that probability limits do
not exist for most of the problems previously leading to
inconsistencies.  Thus, the MP remains a serious challenge to
objective Bayes.

\section{The Marginalization Paradox}
\label{sec:mp}

We briefly describe the key elements of the MP.  Let $p(x|\theta)$ be
the density function for some statistical model, $\pi(\theta)$ a
possibly improper prior, and $\pi(\theta|x)$ the corresponding
posterior, as computed formally from Bayes' law. As noted above, if
$\pi(\theta)$ is improper, i.e., if its integral is infinite, we call
$\pi(\theta|x)$ a \textit{formal} posterior.  Now suppose that
$x=(y,z)$ and $\theta=(\eta,\zeta)$, that the marginal density
$p(z|\theta)$ depends on $\theta$ only through $\zeta$, and that the
marginal posterior $\pi(\zeta|x)$ depends on $x$ only through $z$. We
denote these functions by $\ptilde$ and $\pitilde$, respectively.

Intuition would now suggest that
\begin{displaymath}
  \pitilde \propto \ptilde\,\pi(\zeta),
\end{displaymath}
for some function $\pi(\zeta)$. That is, the marginalized quantities
should satisfy a (possibly formal) Bayes' law. In general, however, we
find that they do not. The problems in which the inconsistencies arise
are ordinary problems of statistical inference, although they need
to satisfy certain symmetry properties. However, they are not
pathological or contrived in any way.

\begin{figure}
  \label{fig:scheme}
  \resizebox{18pc}{!}{\includegraphics{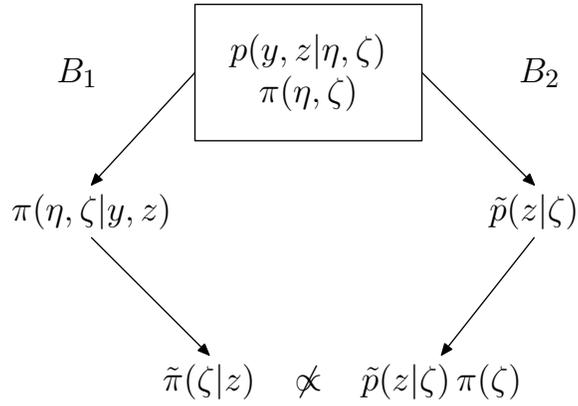}}
\caption{Schematic representation of marginalization paradox.}
\end{figure}

A schematic of the MP is presented in Figure~\ref{fig:scheme}. 
The paradox is sometimes dramatized by ascribing the different
computations to two Bayesians, $B_1$ and $B_2$; the routes they take
are indicated in the Figure. For numerous examples, we refer the
reader to~\cite{SD72,DSZ73}.

\section{Significance of the MP}
\label{sec:ob}

Before proceeding to our analysis, we discuss in greater detail the
significance of the MP for objective Bayes.  One of the perceived
weaknesses in the Bayesian position has been the subjectivity inherent
in using ``personal opinion'' to define the prior distribution.
Objective Bayes addresses this issue by maintaining that prior
distributions can be objectively associated with specific states of
information. Inferences are only subjective in the sense that
different subjects have different states of information.

Objective Bayes maintains, in particular, that there is a unique
numerical representation of ``complete ignorance.''  The operational
definition of ``complete ignorance'' is that our prior beliefs are
invariant under certain transformations of the parameter space. For
example, if we are ignorant of the scale of the variance, our belief
that $\sigma$ lies in an interval $\Delta\sigma$ is the same as our
belief that it lies in $10\Delta\sigma$. These transformations form a
mathematical group.  It is here that the connection with the MP
becomes apparent. If the symmetry group is noncompact, then any
group-invariant (\textit{Haar}) measure is improper.

Many of the symmetries of interest in statistics are described by
noncompact groups. For example, translational symmetry on $R$, the
real numbers, is expressed by the noncompact translation group, whose
invariant measure is Lebesgue measure, with infinite total measure.
Similarly, scale invariance is described by the noncompact group of
positive reals. The invariant measure is $d\sigma/\sigma$, which
diverges. Multivariate analysis typically involves the general linear
group, which is also noncompact.

Ignorance priors are often used directly, but they also form the basis
for the method of maximum entropy. Here, entropy is really ``relative
entropy,'' which is defined relative to some base measure $\pi_0$:
\begin{displaymath}
      H(\pi,\pi_0) = - \int \pi(\theta) 
      \log\left[{\pi(\theta)\over\pi_0(\theta)}\right]\,d\theta.
\end{displaymath}
In the absence of constraints representing additional information,
maximum entropy is achieved for $\pi(\theta)=\pi_0(\theta)$, which
must thus represent a state of complete ignorance.  If we cannot
define a meaningful ignorance prior, we cannot regard maximum entropy
as an ``objective'' procedure.

Jaynes' clearly recognized the challenge the MP posed for his
approach, and opposed it vociferously. A statement of his views can be
found in Chapter~15 of~\cite{Jaynes03}.  ``[The
MP]$\ldots$seemed to threaten the consistency of all probability
theory.''  ``It has been able to do far more damage to the cause of
scientific inference than any other [paradox].''  ``Scientific
inference thus suffered a setback from which it will take decades to
recover.''

Jaynes claims to identify several errors in the analysis of DSZ, but
his basic criticism is that the analysis uses, at critical junctures,
intuitive reasoning about improper quantities, and that the correct
result can only be obtained by carefully considering limits of proper
quantities, although this analysis does not appear to have been
carried out.  He claims that a correct analysis will reveal that $B_1$
and $B_2$ have used different prior information, so it is only to be
expected that their answers will differ.

We will not attempt to address the specific arguments contained
in~\cite{Jaynes03} and elsewhere.  Instead, we start from a point of
agreement between Jaynes and DSZ, and show that the whole
interpretation of the paradox hinges on how the specifics of that
point are interpreted.

\section{Limit concepts}
\label{sec:lc}

The point of agreement is this: both Jaynes and DSZ agree that
``infinite'' quantities must be interpreted as limits of finite
quantities. This is the main motif of Jaynes' Chapter 15: ``The
paradoxes of probability theory.'' It is expressed, for example, in
the statement that ``An improper pdf has meaning only as the limit of
a well-defined sequence of proper pdfs.''~\cite[p.~487]{Jaynes03}.
DSZ are in complete agreement, and say so explicitly: ```Infinity'
finds practical justification only when it can be interpreted as an
idealized approximation of the finite.''~\cite[p.~4]{DSZ96}

In particular, both Jaynes and DSZ agree that if $\pi(\theta)$ is an
improper prior, then in order to define the corresponding posterior,
$\pi(\theta|x)$, we need to construct a sequence of proper priors,
$\{\pi_n(\theta)\}$, such that $\pi_n(\theta)\rightarrow\pi(\theta)$
in some sense, and define $\pi(\theta|x)$ as the limit
\begin{displaymath}
  \pi(\theta|x) \equiv \lim_n \pi_n(\theta|x),
\end{displaymath}
where each $\pi_n(\theta|x)$ is the ordinary posterior
corresponding to $\pi(\theta|x)$. There is, however, more
than one way to define a limit! 

Jaynes follows Jeffreys in adopting what may be called pointwise
limits. We say that $\pipt(\theta|x)$ is a \textit{pointwise limit}
for $p(x|\theta)$ and $\{\pi_n(\theta)\}$ if, \textit{for each $x$},
\begin{equation}
   \label{eq:ptwise}
 \int|\pi_n(\theta|x) - \pipt(\theta|x)|\,d\theta \rightarrow 0.    
\end{equation}
Jaynes adopts this definition explicitly on p.~471 of~\cite{Jaynes03}.
(He does not specify the sense in which the measures
$\pi_n(\theta|x)\,d\theta$ must converge to the measure
$\pipt(\theta|x)\,d\theta$, although this is inessential; in
\eqref{eq:ptwise}, we have used convergence in total variation norm.)
Jeffreys also adopts this definition, although always
implicitly. Thus, for example, Jeffreys writes that ``If \textit{in an
  actual series of observations} the standard deviation is much more
than the smallest admissible value of $\sigma$, and much less than the
largest, the truncation of the distribution makes a negligible change
in the results''~\cite[p.~121]{Jeffreys61}~(italics mine).  A similar
example is worked out explicitly in~\cite[p.~68]{Jeffreys57}, where
Jeffreys calculates a posterior as the limit of $\pi_n(\theta|x)$ with
$x$ fixed (in our notation).

Pointwise limits lead to the formal posterior, and both Jaynes and
Jeffreys justify their use of formal posteriors on this
basis. Typically this is done for specific examples, but in fact, the
argument can be used to justify the formal posterior for essentially
any reasonable prior.  Cf.~\cite{Wallace59} for a general proof under
weak assumptions.

An alternative notion of limit is due to
M.~Stone~\cite{Stone63,Stone65,Stone70}.
%cf.~also~\cite{Stein65,Akaike80} for similar notions. 
We say that
$\pipr(\theta|x)$ is a \textit{probability limit} for $p(x|\theta)$
and $\{\pi_n(\theta)\}$ if
\begin{equation}
  \label{eq:pl}
  \int\left[\int|\pi_n(\theta|x) - \pipr(\theta|x)|\,d\theta\right] 
  \,m_n(x) \,dx\rightarrow 0,
\end{equation}
where $m_n(x)$ is the marginal data density, $\int
p(x|\theta)\,\pi_n(\theta)\,d\theta$. 
Thus, pointwise limits require that, for any fixed $x$, the bracketed
expression eventually becomes small; probability limits require that
the \textit{average} of this expression over $m_n(x)$ eventually
becomes small. The intuition behind this definition is that the true
prior is close to one of the $\pi_n(\theta)$, and $\pipr(\theta|x)$ is
an idealization. It is a useful idealization if it would usually be a
good approximation in the region where the data is expected.

At first glance, both definitions seem reasonable, and it may seem
that the difference could at most be technical. In fact, the
difference is profound, as will become apparent
when we examine a particular example.

\section{Stone's example}
\label{sec:se}

In this section, I present an example, due to Stone~\cite{Stone76},
which illustrates the difference between the two limit concepts, and
also illustrates some disturbing properties of the pointwise
concept. Consider a Gaussian random variable, $X\sim N(\theta,1)$, and
assume that, rather than using the usual uniform prior on $\theta$, we
use an exponential: $\pi(\theta) = \exp(a\theta)$. Although no one
would use this prior in practice, it illustrates phenomena that arise
with realistic priors in more complicated problems~\cite{Stone76}.

We may approximate the improper exponential prior with the sequence
$\pi_n(\theta)$:
\begin{equation}
  \label{eq:prop}
      \pi_n(\theta) \propto
      \exp\left[- {(\theta - an)^2\over 2n}\right]
      \propto  \exp\left[ a\theta - {\theta^2\over 2n} \right].  
\end{equation}
Pointwise limits and probability limits both exist for this
sequence, but they are different:
\begin{equation}
  \label{eq:ptvspr}
  \pipt(\theta|x) \propto
      \exp[-\onehalf{(\theta - x - a)^2}], \qquad
\pipr(\theta|x) \propto
      \exp[-\onehalf{(\theta - x)^2}].
\end{equation}
The first expression is just the formal posterior; the second
result follows from Stone's theorem~\cite{Stone70}.

\begin{figure}
  \label{fig:plot}
  \resizebox{20pc}{!}{\includegraphics{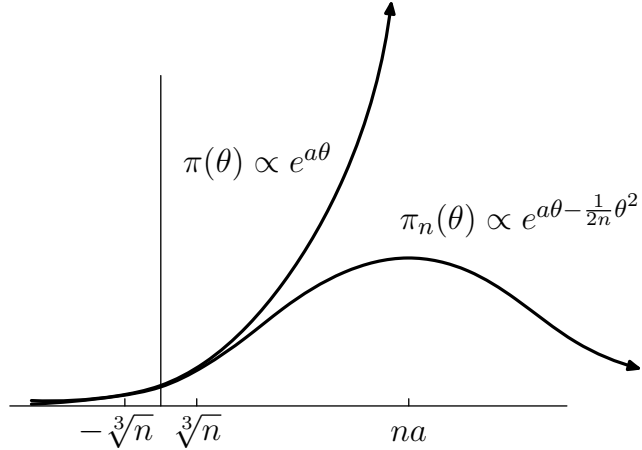}}
\caption{Comparison of improper prior and approximation.}
\end{figure}

To understand the difference between pointwise and probability limits,
consider Figure~\ref{fig:plot}, 
which compares the priors $\pi(\theta)$ and $\pi_n(\theta)$. The
likelihood is local, so the posteriors $\pi(\theta|x)$ and
$\pi_n(\theta|x)$ will be similar if the priors are nearly
proportional in the vicinity of $x$. Consider the interval $I_n =
(-\root 3 \of n, \root 3 \of n)$. As $n\rightarrow\infty$ the priors,
which differ only by the term $\theta^2/n= O(1/{\root 3 \of n})$ in
the exponential, will converge in $I_n$, and $\pi_n(\theta|x)$ will
converge to $\pi(\theta|x)$ for each $x$ in $I_n$. As
$n\rightarrow\infty$, however, the interval $I_n$ will expand to cover
the entire real line. Thus, the pointwise limit of the
$\pi_n(\theta|x)$, which is the formal posterior, will exist. 

But is this a good reason to regard the formal posterior as the limit
of the $\pi_n(\theta|x)$? Consider this from the standpoint of someone
whose true prior is $\pi_n(\theta)$, and who seeks an idealized
posterior that he can use as a good approximation. As we have just seen,
the assertion that $\pi_n(\theta|x)\rightarrow\pipt(\theta|x)$ is
based entirely on agreement in the region close to the origin, where
the data will almost never be found! On the other hand, in the region
where the data \textit{is} expected, $x\sim na$, $\pi_n(\theta|x)$
differs markedly from $\pipt(\theta|x)$. (Explicit formulas are given
in the next section.) In this sense, the formal posterior is poor
approximation to $\pi_n(\theta|x)$, for any $n$.

In the next section, we show that $\pipr(\theta|x)$, which does not
obey the formal Bayes' law, \textit{is} a good approximation in
the region where the data expected.

\section{Bayes' law: local but not global}

We infer from \eqref{eq:pl} that  where $m_n(x)$ is large,
  \begin{equation}
    \label{eq:local}
  \pipr(\theta|x)\approx \pi_n(\theta|x) \propto p(\theta|x)\,\pi_n(\theta).  \end{equation}
That is, $\pipr(\theta|x)$ satisfies Bayes' law locally, and
approximately. On the other hand, in regions where $m_n(x)$ is small,
there is no need for $\pipr(\theta|x)$ to satisfy Bayes' law, and in
general, it will not.

We examine this phenomenon using Stone's example.  First, we show that
Bayes' law holds locally for $\pipr(\theta|x)$, as given in
Eq.~\eqref{eq:ptvspr}; this justifies our assertion that $\pipr(\theta|x)$ is
a probability limit for this problem. It is sufficient to show that
$\pipr(\theta|x)\approx \pi_n(\theta|x)$ when $m_n(x)$ is large.  The
densities $\pi_n(\theta|x)$ and $m_n(x)$ are easily calculated. Let
$\sigma_n^2 = n/(1+n)$. Then
\begin{displaymath}
    \pi_n(\theta|x) \propto 
    \exp\left[-{\left(\theta - \sigma_n^2(x+a)\right)^2
        \over 2\sigma_n^2}\right];\qquad  
    m_n(x) \propto \exp\left[-{(x-an)^2\over 2(1+n)}\right].
\end{displaymath}
Let $x = an+\epsilon$. Then
\begin{displaymath}
    \pi_n(\theta|x) \propto 
    \exp\left[-{\left(\theta - x + {\epsilon\over n+1}\right)^2
        \over 2\sigma_n^2}\right].
\end{displaymath}
In the region where $m_n(x)$ is large, $\epsilon= O(\sqrt{n+1})$.  As
$n\rightarrow\infty$, $\epsilon/(n+1)\rightarrow0$ and
$\sigma_n\rightarrow1$, so $\pipr(\theta|x)\approx \pi_n(\theta|x)$,
as was to be shown.

It is obvious, however, that Bayes' law does not hold globally for
$\pipr(\theta|x)$:
\begin{eqnarray*}
  e^{-\onehalf(x-\theta)^2} &\not\propto&   e^{-\onehalf(x-\theta)^2} e^{a\theta}
\end{eqnarray*}

\section{Discussion}
\label{sec:disc}

The presence or absence of a MP is determined by the type of limit we
use to define improper inferences.  In~\cite{Wallstrom07}, it is shown
that if $\pi(\theta|x)$ is a probability limit, then the marginal,
$\pi(\zeta|z)$, is also a probability limit. The MP shows that if
$\pi(\theta|x)$ is a formal posterior, then the marginal,
$\pi(\zeta|z)$, is not necessarily a formal posterior.  Thus, the
interpretation of improper inferences as probability limits is
internally consistent, at least with respect to the MP, whereas their
interpretation as formal posteriors is not.

The use of probability limits, rather than pointwise limits, is
already known to resolve other difficulties, such as ``strong
inconsistency''~\cite{Stone76}, and ``incoherence.''~\cite{EF04}. Now
that probability limits have been shown to resolve the MP as well, we
have a common explanation for most of the key difficulties of improper
inference: they arise because the pointwise limit, and the formal
posteriors to which they give rise, are fundamentally
unsound. Probability limits, by contrast, appear to provide a
consistent theory of improper inference.

The local nature of Bayes' law for probability limits is the key to
understanding the MP. The requirement that $\pi(\theta|x)$ be a
probability limit is essentially equivalent to the requirement that
Bayes' law hold locally, where $m_n(x)$ is large, in the sense of
\eqref{eq:local}. The requirement that $\pi(\theta|x)$ be a formal
posterior is equivalent to the requirement that Bayes' law hold
globally. $B_1$ and $B_2$ both make inferences based on the latter
requirement. Thus, each of them make inferences that are sometimes
erroneous, and so they frequently disagree.

In resolving the MP, we have followed Jaynes in regarding an improper
inference as a limit of ordinary inferences, based on proper
priors. There is nothing in our procedure that is inconsistent with
the rules of probability theory, as developed by Jaynes in Chapter 2
of~\cite{Jaynes03}. These rules, however, do not apply to the case
where the prior is improper, nor do they stipulate how the improper
case is to be regarded as a limit. In particular, they do not justify
the use of formal posteriors. The use of pointwise limits is an
additional assumption, the validity of which is open to
question. Following Stone, we have argued that pointwise limits are
unsound, and that probability limits better capture the
intuitive meaning of convergence.

It remains to discuss the implications of this analysis for objective
Bayesianism. There are strong indications that the requirement that an
improper inference be a probability limit is very restrictive.  In
group models, Stone has shown that the formal posterior can only be a
probability limit if the prior is right Haar measure and the group
satisfies a technical condition, known as amenability~\cite{Stone70}.
Eaton and Sudderth have shown that many of the formal posteriors of
multivariate analysis are ``incoherent'' or strongly
inconsistent, and thus cannot be probability limits~\cite{ES93}.

Probability limits may be used to construct improper inferences for
the translation and scale groups, and these coincide with the formal
posteriors. Since the most common applications of improper inference
involve these groups, our analysis shows why formal posteriors appear
to work in these simple situations. When we get to more complicated
problems, however, such as those of multivariate analysis, it appears
that probability limits associated with the ``ignorance priors'' for
the relevant symmetries, such as $GL(n)$, do not exist.  Thus, the use
of probability limits restores the technical viability of objective
Bayes to only a very limited domain of improper problems.

\begin{theacknowledgments}
I acknowledge support from
the Department of Energy under contract DE-AC52-06NA25396.
\end{theacknowledgments}

\end{document}